\documentclass[10pt]{article}
\usepackage{amstext}
\usepackage[fleqn]{amsmath}
\usepackage{amsfonts}
\usepackage{amssymb}
\usepackage{amsthm}
\usepackage{newlfont}
\usepackage{graphicx}
\usepackage{sectsty}
\theoremstyle{plain} \theoremstyle{theorem}
\newtheorem{theorem}{Theorem}[section]
\theoremstyle{example}

\theoremstyle{corollary}
\newtheorem{corollary}{Corollary}[section]
\theoremstyle{lemma}
\newtheorem{lemma}{Lemma}[section]
\theoremstyle{proposition}

\theoremstyle{axiom}

\theoremstyle{notation}

\theoremstyle{fact}

\theoremstyle{definition}
\newtheorem{definition}{Definition}[section]
\theoremstyle{remark}
\newtheorem{remark}{Remark}[section]
\numberwithin{equation}{section}
\begin{document}
\title{\bf \Large  An extension of basic Humbert hypergeometric functions $\mathbf{\Phi}_{1}$, $\mathbf{\Phi}_{2}$ and $\mathbf{\Phi}_{3}$}
\author{Ayman Shehata \thanks{%
E-mail: aymanshehata@science.aun.edu.eg,\;drshehata2006@yahoo.com,\; drshehata2009@gmail.com}\\
{\small Department of Mathematics, Faculty of Science, Assiut University, Assiut 71516, Egypt.}}
\date{}
\maketitle{}
\begin{abstract}
Given the growing quantity of proposals and works of basic hypergeometric functions in the scope of $q$-calculus, it is important to introduce a systematic classification of $q$-calculus. Our aim in this article is to investigate several interesting $q$-partial derivative formulas, $q$-contiguous function relations, $q$-recurrence relations, various $q$-partial differential equations, summation formulas, transformation formulas and $q$-integrals representations for basic Humbert hypergeometric functions $\mathbf{\Phi}_{1}$, $\mathbf{\Phi}_{2}$ and $\mathbf{\Phi}_{3}$ under constraints of symmetry parameters. These interesting properties, as special cases, include many known expansions of basic Humbert hypergeometric functions $\mathbf{\Phi}_{1}$, $\mathbf{\Phi}_{2}$ and $\mathbf{\Phi}_{3}$, and are also include particular interest in the area.
\end{abstract}
\section{Introduction and Definitions}
Basic hypergeometric series have assumed great importance during the past decades or so due to their applications in diverse areas. That means, for example, number theory, partition theory,  statistics, combinatorial analysis, quantum mechanics and vector spaces, etc~\cite{ag1, e2, gr1, ka, sr1, rc, sa, ve1}. The development of the transformation theory  is based on the summation formulas for basic hypergeometric series. Indeed, most of the summation formulas for basic hypergeometric series are proved as special cases of the transformation of the hypergeometric series, which are used again to obtain new transformation formulas and to prove the various summation formulas for basic hypergeometric series in a scientific manner, without the use of the majority of the transformation theory of basic hypergeometric series (see for example Jain~\cite{ja1, ja2}, Kandu~\cite{ka}, Rahman and Verma~\cite{rv}, Upadhyay~\cite{up}, Verma~\cite{ve1}, Verma and Jain~\cite{vj1, vj2}, Wang and Chern~\cite{wc}). Later on Jackson~\cite{j3, j4, j5}, Agarwal~\cite{ag1}  defined and introduced the $q$-analogue of Appell functions. In~\cite{krl}, Kim et~al. studied the $q$-analogue of Kummer’s theorem and its contiguous results. In~\cite{hkrr,hkrp}, Harsh et al. investigated of $q$-contiguous functional relations of basic hypergeometric functions. In~\cite{an1}, Andrews derived the summations and transformations for basic Appell series. In~\cite{e1, e2, e3, e4, e5, e6, e7, e8}, Ernst discussed the certain generalizations of $q$-hypergeometric, $q$-Appell and $q$-Lauricella series of two variables. In~\cite{ja1, ja2}, Jain established the transformations of basic hypergeometric series and their applications. In~\cite{ve1, vj1, vj2}, Verma  and Jain introduced the transformations of basic hypergeometric series. In~\cite{sh1, sh2}, Shehata presented the $(p,q)$-Bessel and $(p,q)$-Humbert functions. Several properties for basic Horn functions $H_{3}$, $H_{4}$, $H_{6}$ and $H_{7}$ was earlier investigated by Shehata~\cite{sh3, sh4}. Recently, AL E'damat and Shehata \cite{as}, and Qing-Bo Cai et al. \cite{qgsa} investigated some properties of bibasic Humbert hypergeometric functions $\Phi_{1}$, $\Xi_{1}$ and $\Xi_{2}$. Shehata et al. \cite{ssps} have investigated and discussed the new formulas for basic Humbert hypergeometric functions.
\begin{definition}
Let $0<|q|<1$, $q\in\mathbb{C}$ and $\ell\geq 0$, the $q$-shifted factorial is defined by (see~\cite{sh3, sh4})
\begin{eqnarray}
\begin{split}
&(q^{\alpha};q)_{\ell}=\left\{
            \begin{array}{ll}
              \prod_{r=0}^{\ell-1}(1- q^{\alpha+r}), & \hbox{$\ell\geq1,q^{\alpha}\in\mathbb{C}\backslash\{ 1, q^{-1}, q^{-2},\ldots,q^{1-\ell}\}$;} \\
              1, & \hbox{$\ell=0,q^{\alpha}\in\mathbb{C}$.}
            \end{array}
          \right.\\
          &=\left\{
  \begin{array}{ll}
    (1-q^{\alpha})(1-q^{\alpha+1})\ldots(1-q^{\alpha+\ell-1}), & \hbox{$\ell\in\Bbb{N},q^{\alpha}\in\mathbb{C}\backslash\{ 1, q^{-1}, q^{-2},\ldots,q^{1-\ell}\}$;} \\
    1, & \hbox{$\ell=0,q^{\alpha}\in\mathbb{C}$,}
  \end{array}
\right.\label{1.1}
\end{split}
\end{eqnarray}
where the symbols $\mathbb{C}$, $\Bbb{N}$ denote the set of complex numbers and natural integers, respectively.
\end{definition}
\begin{definition}
For $Re(m)>0$, $n\neq 0,-1,-2,\ldots$, the $q$-beta function is given as follows (see~\cite{gr1}) 
\begin{eqnarray}
\begin{split}
B_{q}(m,n)=\int_{0}^{1}t^{m-1}\frac{(qt;q)_{\infty}}{(tq^{n};q)_{\infty}}d_{q}t.\label{1.2}
\end{split}
\end{eqnarray}
\end{definition}  
The notations and identities which we need in our subsequent work are as follows: (see~\cite{sw})
\begin{eqnarray}
\begin{split}
(q^{\alpha+1};q)_{\ell+\kappa}&=\frac{1-q^{\alpha+\ell+\kappa}}{1-q^{\alpha}}(q^{\alpha};q)_{\ell+\kappa}=\frac{[\alpha+\ell+\kappa]_{q}}{[\alpha]_{q}}(q^{\alpha};q)_{\ell+\kappa},\ell,\kappa\in\Bbb{N}_{0}=\Bbb{N}\cup\{0\},\label{1.3}
\end{split}
\end{eqnarray}
where the $q$-number (basic or quantun number) $[\alpha]_{q}$ is defined by (see~\cite{gr1})
\begin{eqnarray*}
\begin{split}
[\alpha]_{q}=\frac{1-q^{\alpha}}{1-q}, \alpha\in \mathbb{C}, 0<|q|<1, q\in\mathbb{C},
\end{split}
\end{eqnarray*}\vspace{-12pt}
\begin{eqnarray}
\begin{split}
(q^{\alpha+1};q)_{\ell+k}&=\bigg{[}1+q^{\alpha}\frac{1-q^{\ell}}{1-q^{\alpha}}+q^{\alpha+\ell}\frac{1-q^{k}}{1-q^{\alpha}}\bigg{]}(q^{\alpha};q)_{\ell+k},\\
&=\bigg{[}1+q^{\alpha}\frac{1-q^{k}}{1-q^{\alpha}}+q^{\alpha+k}\frac{1-q^{\ell}}{1-q^{\alpha}}\bigg{]}(q^{\alpha};q)_{\ell+k},q^{\alpha}\neq 1,\ell,k\in\Bbb{N}_{0}\label{1.4}
\end{split}
\end{eqnarray}\vspace{-12pt}
and
\begin{eqnarray}
\begin{split}
\frac{1}{(q^{\alpha-1};q)_{\ell+\kappa}}&=\bigg{(}1+q^{\alpha}\frac{1-q^{\kappa}}{q-q^{\alpha}}+q^{\alpha+\kappa}\frac{1-q^{\ell}}{q-q^{\alpha}}\bigg{)}\frac{1}{(q^{\alpha};q)_{\ell+\kappa}}\\
&=\bigg{(}1+q^{\alpha}\frac{1-q^{\ell}}{q-q^{\alpha}}+q^{\alpha+\ell}\frac{1-q^{\kappa}}{q-q^{\alpha}}\bigg{)}\frac{1}{(q^{\alpha};q)_{\ell+\kappa}},q^{\alpha}\neq q,\ell,\kappa\in\Bbb{N}_{0}.\label{1.5}
\end{split}
\end{eqnarray}
\begin{definition}
The $q$-difference operator $D_{x,q}$ is defined by Jackson~\cite{j3, j4, j5} as follows
\begin{equation}
\begin{split}
D_{\mathrm{z},q}\Psi(\mathrm{z})=\frac{\Psi(\mathrm{z})-\Psi(q\mathrm{z})}{(1-q)\mathrm{z}},\mathrm{z}\neq0,\label{1.6}
\end{split}
\end{equation}
where $\Psi$ is a function defined on a subset of the complex or real plane.
\end{definition}
\begin{definition}
For $0<|q|<1$ and $q\in\mathbb{C}$, we define the basic Humbert hypergeometric functions $\mathbf{\Phi}_{1}$, $\mathbf{\Phi}_{2}$ and $\mathbf{\Phi}_{3}$ as follows in the form
\begin{eqnarray}
\begin{split}
&\mathbf{\Phi}_{1}=\mathbf{\Phi}_{1}(q^{\alpha},q^{\beta};q^{\gamma};q,x,y)=\sum_{n,k=0}^{\infty}\frac{(q^{\alpha};q)_{n+k}(q^{\beta};q)_{n}x^{n}y^{k}}{(q^{\gamma};q)_{n+k}(q;q)_{n}(q;q)_{k}},q^{\gamma}\neq 1, q^{-1}, q^{-2},\ldots,\label{1.7}
\end{split}
\end{eqnarray}\vspace{-12pt}
\begin{eqnarray}
\begin{split}
&\mathbf{\Phi}_{2}=\mathbf{\Phi}_{2}(q^{\alpha},q^{\beta};q^{\gamma};q,x,y)=\sum_{n,k=0}^{\infty}\frac{(q^{\alpha};q)_{n}(q^{\beta};q)_{k}x^{n}y^{k}}{(q^{\gamma};q)_{n+k}(q;q)_{n}(q;q)_{k}},q^{\gamma}\neq 1, q^{-1}, q^{-2},\ldots\label{1.8}
\end{split}
\end{eqnarray}
and
\begin{eqnarray}
\begin{split}
&\mathbf{\Phi}_{3}=\mathbf{\Phi}_{3}(q^{\alpha};q^{\beta};q,x,y)=\sum_{n,k=0}^{\infty}\frac{(q^{\alpha};q)_{n}x^{n}y^{k}}{(q^{\beta};q)_{n+k}(q;q)_{n}(q;q)_{k}},q^{\beta}\neq 1, q^{-1}, q^{-2},\ldots.\label{1.9}
\end{split}
\end{eqnarray}
\end{definition}
\begin{lemma}
For $\ell\geq 0$ and $\kappa\geq 0$, the relation is given by (see~\cite{ra})
\begin{eqnarray}
\begin{split}
\sum_{\ell=0}^{\infty}\sum_{\kappa=0}^{\infty}\mathrm{A}(\ell,\kappa)=\sum_{\ell=0}^{\infty}\sum_{\kappa=0}^{\ell}\mathrm{A}(\ell-\kappa,\kappa).\label{1.10}
\end{split}
\end{eqnarray}
\end{lemma}
In this study we express a family of extended forms, to motivate the former works~\cite{sh3, sh4}. In this paper, in a systematic way, we aim at establishing new $q$-partial derivative formulas, $q$-differential formulas, $q$-recursion formulas, $q$-differential recursion relations, $q$-contiguous function relations, $q$-partial differential equations, summation formulas, transformation formulas, integral representations and interesting transformations for $q$ or basic Humbert hypergeometric functions $\mathbf{\Phi}_{1}$, $\mathbf{\Phi}_{2}$ and $\mathbf{\Phi}_{3}$ under what constraints of symmetry parameters..
\section{Derivation of Main Results}
In this section, we will derive our main results.
\begin{theorem} The following $q$-derivative formulas hold true for $\mathbf{\Phi}_{1}$, $\mathbf{\Phi}_{2}$ and $\mathbf{\Phi}_{3}$
\begin{eqnarray}
\begin{split}
&\mathbb{D}_{x,q}^{r}\mathbf{\Phi}_{1}=\frac{(q^{\alpha};q)_{r}(q^{\beta};q)_{r}}{(q^{\gamma};q)_{r}(1-q)^{r}}\mathbf{\Phi}_{1}(q^{\alpha+r},q^{\beta+r};q^{\gamma+r}),\label{2.1}
\end{split}
\end{eqnarray}\vspace{-12pt}
\begin{eqnarray}
\begin{split}
&\mathbb{D}_{y,q}^{s}\mathbf{\Phi}_{1}=\frac{(q^{\alpha};q)_{s}}{(q^{\gamma};q)_{s}(1-q)^{s}}\mathbf{\Phi}_{1}(q^{\alpha+s};q^{\gamma+s}),\\
&\mathbb{D}_{x,q}^{r}\mathbb{D}_{y,q}^{s}\mathbf{\Phi}_{1}=\frac{(q^{\alpha};q)_{r+s}(q^{\beta};q)_{r}}{(q^{\gamma};q)_{r+s}(1-q)^{r+s}}\mathbf{\Phi}_{1}(q^{\alpha+r+s},q^{\beta+r};q^{\gamma+r+s}),\label{2.2}
\label{2.10}
\end{split}
\end{eqnarray}\vspace{-12pt}
\begin{eqnarray}
\begin{split}
&\mathbb{D}_{x,q}^{r}\mathbf{\Phi}_{2}=\frac{(q^{\alpha};q)_{r}}{(q^{\gamma};q)_{r}(1-q)^{r}}\mathbf{\Phi}_{2}(q^{\alpha+r};q^{\gamma+r}),\\
&\mathbb{D}_{y,q}^{s}\mathbf{\Phi}_{2}=\frac{(q^{\beta};q)_{s}}{(q^{\gamma};q)_{s}(1-q)^{s}}\mathbf{\Phi}_{2}(q^{\beta+s};q^{\gamma+s}),\\
&\mathbb{D}_{x,q}^{r}\mathbb{D}_{y,q}^{s}\mathbf{\Phi}_{2}=\frac{(q^{\alpha};q)_{r}(q^{\beta};q)_{s}}{(q^{\gamma};q)_{r+s}(1-q)^{r+s}}\mathbf{\Phi}_{2}(q^{\alpha+r},q^{\beta+s};q^{\gamma+r+s})\label{2.3}
\label{3.10}
\end{split}
\end{eqnarray}\vspace{-12pt}
and
\begin{eqnarray}
\begin{split}
&\mathbb{D}_{x,q}^{r}\mathbf{\Phi}_{3}=\frac{(q^{\alpha};q)_{r}}{(q^{\beta};q)_{r}(1-q)^{r}}\mathbf{\Phi}_{3}(q^{\alpha+r};q^{\beta+r}),\\
&\mathbb{D}_{y,q}^{s}\mathbf{\Phi}_{3}=\frac{1}{(q^{\beta};q)_{s}(1-q)^{s}}\mathbf{\Phi}_{3}(q^{\beta+s}),\\
&\mathbb{D}_{x,q}^{r}\mathbb{D}_{y,q}^{s}\mathbf{\Phi}_{3}=\frac{(q^{\alpha};q)_{r}}{(1-q)^{r+s}(q^{\beta};q)_{r+s}}\mathbf{\Phi}_{3}(q^{\alpha+r};q^{\beta+r+s}).\label{2.4}
\label{4.10}
\end{split}
\end{eqnarray}
\end{theorem}
\begin{proof}
Using (\ref{1.6}) and (\ref{1.7}), by differentiating of $\mathbf{\Phi}_{1}$ with respect to $x$, we obtain
\begin{eqnarray}
\begin{split}
&\mathbb{D}_{x,q}\mathbf{\Phi}_{1}=\sum_{n,k=0}^{\infty}\frac{1-q^{n}}{1-q}\frac{(q^{\alpha};q)_{n+k}(q^{\beta};q)_{n}x^{n}y^{k}}{(q^{\gamma};q)_{n+k}(q;q)_{k}(q;q)_{n}}=\frac{(1-q^{\alpha})(1-q^{\beta})}{(1-q)(1-q^{\gamma})}\mathbf{\Phi}_{1}(q^{\alpha+1},q^{\beta+1};q^{\gamma+1}),q^{\gamma}\neq 1.\label{2.5}
\end{split}
\end{eqnarray}

Repeating the $q$-derivatives of $\mathbf{\Phi}_{1}$ with respect to $x$ for $r$ times, we obtain (\ref{2.1}). In a similar way, we obtain (\ref{2.2})--(\ref{2.4}) .
\end{proof}

\begin{corollary} The following relations for $\mathbf{\Phi}_{1}$, $\mathbf{\Phi}_{2}$ and $\mathbf{\Phi}_{3}$ hold true
\begin{eqnarray}
\begin{split}
q^{\mathbf{\Theta}_{x,q}}\mathbf{\Phi}_{1}(q^{\alpha},q^{\beta};q^{\gamma};q,x,y)=\mathbf{\Phi}_{1}(xq),\\
q^{\mathbf{\Theta}_{y,q}}\mathbf{\Phi}_{1}(q^{\alpha},q^{\beta};q^{\gamma};q,x,y)=\mathbf{\Phi}_{1}(yq),\label{2.6}
\end{split}
\end{eqnarray}\vspace{-12pt}
\begin{eqnarray}
\begin{split}
q^{\mathbf{\Theta}_{x,q}}\mathbf{\Phi}_{2}(q^{\alpha},q^{\beta};q^{\gamma};q,x,y)=\mathbf{\Phi}_{2}(xq),\\
q^{\mathbf{\Theta}_{y,q}}\mathbf{\Phi}_{2}(q^{\alpha},q^{\beta};q^{\gamma};q,x,y)=\mathbf{\Phi}_{2}(yq)\label{2.7}
\end{split}
\end{eqnarray}\vspace{-12pt}
and
\begin{eqnarray}
\begin{split}
q^{\mathbf{\Theta}_{x,q}}\mathbf{\Phi}_{3}(q^{\alpha},q^{\beta};q^{\gamma};q,x,y)=\mathbf{\Phi}_{3}(xq),\\
q^{\mathbf{\Theta}_{y,q}}\mathbf{\Phi}_{3}(q^{\alpha},q^{\beta};q^{\gamma};q,x,y)=\mathbf{\Phi}_{3}(yq),\label{2.8}
\end{split}
\end{eqnarray}
where $\mathbf{\Theta}_{x,q}=x\mathbb{D}_{x,q}$ and $\mathbf{\Theta}_{y,q}=y\mathbb{D}_{y,q}$.
\end{corollary}
\begin{proof}
Using an operator $q^{\mathbf{\Theta}_{x,q}}f(x)=q^{\mathbf{x\mathbb{D}}_{x,q}}f(x)=f(qx)$, we obtain (\ref{2.6})--(\ref{2.8}).
\end{proof}

\begin{theorem} The functions $\mathbf{\Phi}_{1}$, $\mathbf{\Phi}_{2}$  and $\mathbf{\Phi}_{3}$ satisfy the $q$-differential formulas:
\begin{eqnarray}
\begin{split}
&\bigg{(}q^{\alpha}\mathbf{\Theta}_{x,q}+[\alpha]_{q}\bigg{)}\mathbf{\Phi}_{1}+q^{\alpha}\mathbf{\Theta}_{y,q}\mathbf{\Phi}_{1}(qx)=[\alpha]_{q}\mathbf{\Phi}_{1}(q^{\alpha+1}),\label{2.9}
\end{split}
\end{eqnarray}\vspace{-12pt}
\begin{eqnarray}
\begin{split}
&\bigg{(}q^{\alpha}\mathbf{\Theta}_{y,q}+[\alpha]_{q}\bigg{)}\mathbf{\Phi}_{1}+q^{\alpha}\mathbf{\Theta}_{x,q}\mathbf{\Phi}_{1}(qy)=[\alpha]_{q}\mathbf{\Phi}_{1}(q^{\alpha+1}),\\
&\bigg{(}q^{\beta}\mathbf{\Theta}_{x,q}+[\beta]_{q}\bigg{)}\mathbf{\Phi}_{1}=[\beta]_{q}\mathbf{\Phi}_{1}(q^{\beta+1}),\\
&\bigg{(}q^{\gamma-1}\mathbf{\Theta}_{x,q}+[\gamma-1]_{q}\bigg{)}\mathbf{\Phi}_{1}+q^{\gamma-1}\mathbf{\Theta}_{y,q}\mathbf{\Phi}_{1}(qx)=[\gamma-1]_{q}\mathbf{\Phi}_{1}(q^{\gamma-1}),\\
&\bigg{(}q^{\gamma-1}\mathbf{\Theta}_{y,q}+[\gamma-1]_{q}\bigg{)}\mathbf{\Phi}_{1}+q^{\gamma-1}\mathbf{\Theta}_{x,q}\mathbf{\Phi}_{1}(qy)=[\gamma-1]_{q}\mathbf{\Phi}_{1}(q^{\gamma-1}),\label{2.10}
\end{split}
\end{eqnarray}\vspace{-12pt}
\begin{eqnarray}
\begin{split}
&\bigg{(}q^{\alpha}\mathbf{\Theta}_{x,q}+[\alpha]_{q}\bigg{)}\mathbf{\Phi}_{2}=[\alpha]_{q}\mathbf{\Phi}_{2}(q^{\alpha+1}),\\
&\bigg{(}q^{\beta}\mathbf{\Theta}_{y,q}+[\beta]_{q}\bigg{)}\mathbf{\Phi}_{2}=[\beta]_{q}\mathbf{\Phi}_{2}(q^{\beta+1}),\\
&\bigg{[}q^{\gamma-1}\mathbf{\Theta}_{x,q}+[\gamma-1]_{q}\bigg{]}\mathbf{\Phi}_{2}+q^{\gamma-1}\mathbf{\Theta}_{y,q}\mathbf{\Phi}_{2}(qx)=[\gamma-1]_{q}\mathbf{\Phi}_{2}(q^{\gamma-1}),\\
&\bigg{[}q^{\gamma-1}\mathbf{\Theta}_{y,q}+[\gamma-1]_{q}\mathbf{\Phi}_{2}+q^{\gamma-1}\mathbf{\Theta}_{x,q}\mathbf{\Phi}_{2}(qy)=[\gamma-1]_{q}\mathbf{\Phi}_{2}(q^{\gamma-1})\label{2.11}
\end{split}
\end{eqnarray}\vspace{-12pt}
and
\begin{eqnarray}
\begin{split}
&\bigg{(}q^{\alpha}\mathbf{\Theta}_{x,q}+[\alpha]_{q}\bigg{)}\mathbf{\Phi}_{3}=[\alpha]_{q}\mathbf{\Phi}_{3}(q^{\alpha+1}),\\
&\bigg{(}q^{\beta-1}\mathbf{\Theta}_{x,q}+[\beta-1]_{q}\bigg{)}\mathbf{\Phi}_{3}+q^{\beta-1}\mathbf{\Theta}_{y,q}\mathbf{\Phi}_{3}(qx)=[\beta-1]_{q}\mathbf{\Phi}_{3}(q^{\beta-1}),\\
&\bigg{(}q^{\beta-1}\mathbf{\Theta}_{y,q}+[\beta-1]_{q}\bigg{)}\mathbf{\Phi}_{3}+q^{\beta-1}\mathbf{\Theta}_{x,q}\mathbf{\Phi}_{3}(qy)=[\beta-1]_{q}\mathbf{\Phi}_{3}(q^{\beta-1}).\label{2.12}
\end{split}
\end{eqnarray}
\end{theorem}
\begin{proof}
From (\ref{1.4}) and using (\ref{1.7}), we obtain the $q$-contiguous relation:
\begin{eqnarray*}
\begin{split}
&\mathbf{\Phi}_{1}(q^{\alpha+1})=\sum_{n,k=0}^{\infty}\frac{1-q^{\alpha+n+k}}{1-q^{\alpha}}\frac{(q^{\alpha};q)_{n+k}(q^{\beta};q)_{n}x^{n}y^{k}}{(q^{\gamma};q)_{n+k}(q;q)_{k}(q;q)_{n}}\\
=&\sum_{n,k=0}^{\infty}\bigg{(}1+q^{\alpha}\frac{1-q^{n}}{1-q^{\alpha}}+q^{\alpha+n}\frac{1-q^{k}}{1-q^{\alpha}}\bigg{)}\frac{(q^{\alpha};q)_{n+k}(q^{\beta};q)_{n}x^{n}y^{k}}{(q^{\gamma};q)_{n+k}(q;q)_{k}(q;q)_{n}}\\
=&\sum_{n,k=0}^{\infty}\frac{(q^{\alpha};q)_{n+k}(q^{\beta};q)_{n}x^{n}y^{k}}{(q^{\gamma};q)_{n+k}(q;q)_{k}(q;q)_{n}}+\frac{q^{\alpha}(1-q)}{1-q^{\alpha}}\sum_{n,k=0}^{\infty}\frac{1-q^{n}}{1-q}\frac{(q^{\alpha};q)_{n+k}(q^{\beta};q)_{n}x^{n}y^{k}}{(q^{\gamma};q)_{n+k}(q;q)_{k}(q;q)_{n}}\\
+&\frac{q^{\alpha}(1-q)}{1-q^{\alpha}}\sum_{n,k=0}^{\infty}\frac{1-q^{k}}{1-q}\frac{(q^{\alpha};q)_{n+k}(q^{\beta};q)_{n}(qx)^{n}y^{k}}{(q^{\gamma};q)_{n+k}(q;q)_{k}(q;q)_{n}}\\
=&\mathbf{\Phi}_{1}(q^{\alpha},q^{\beta};q^{\gamma};q,x,y)+\frac{q^{\alpha}}{[\alpha]_{q}}\mathbf{\Theta}_{x,q}\mathbf{\Phi}_{1}(q^{\alpha},q^{\beta};q^{\gamma};q,x,y)+\frac{q^{\alpha}}{[\alpha]_{q}}\mathbf{\Theta}_{y,q}\mathbf{\Phi}_{1}(q^{\alpha},q^{\beta};q^{\gamma};q,qx,y).
\end{split}
\end{eqnarray*}

Proofs (\ref{2.10})--(\ref{2.12}) are similar in technique as the proof (\ref{2.9}).
\end{proof}
\begin{corollary}
The $q$-partial derivatives for functions $\mathbf{\Phi}_{1}$, $\mathbf{\Phi}_{2}$ and $\mathbf{\Phi}_{3}$ hold true
\begin{eqnarray}
\begin{split}
&\mathbf{\Theta}_{x}\mathbf{\Phi}_{1}+\mathbf{\Theta}_{y}\mathbf{\Phi}_{1}(qx)=\mathbf{\Theta}_{y}\mathbf{\Phi}_{1}+\mathbf{\Theta}_{x}\mathbf{\Phi}_{1}(qy),\\
&\big{(}q^{\alpha}-q^{\gamma-1}\big{)}\mathbf{\Phi}_{1}=(1-q^{\gamma-1})q^{\alpha}\mathbf{\Phi}_{1}(q^{\gamma-1})-(1-q^{\alpha})q^{\gamma-1}\mathbf{\Phi}_{1}(q^{\alpha+1}),\label{2.13}
\end{split}
\end{eqnarray}\vspace{-12pt}
\begin{eqnarray}
\begin{split}
&\frac{q^{\alpha}}{[\alpha]_{q}}\mathbf{\Theta}_{x,q}\mathbf{\Phi}_{2}+\mathbf{\Phi}_{2}(q^{\beta+1})=\frac{q^{\beta}}{[\beta]_{q}}\mathbf{\Theta}_{y,q}\mathbf{\Phi}_{2}+\mathbf{\Phi}_{2}(q^{\alpha+1}),\\
&\mathbf{\Theta}_{x,q}\mathbf{\Phi}_{2}+\mathbf{\Theta}_{y,q}\mathbf{\Phi}_{2}(qx)=\mathbf{\Theta}_{y,q}\mathbf{\Phi}_{2}+\mathbf{\Theta}_{x,q}\mathbf{\Phi}_{2}(qy)\label{2.14}
\end{split}
\end{eqnarray}\vspace{-12pt}
and
\begin{eqnarray}
\begin{split}
&\mathbf{\Theta}_{x,q}\mathbf{\Phi}_{3}+\mathbf{\Theta}_{y,q}\mathbf{\Phi}_{3}(qx)=\mathbf{\Theta}_{y,q}\mathbf{\Phi}_{3}+\mathbf{\Theta}_{x,q}\mathbf{\Phi}_{3}(qy).\label{2.15}
\end{split}
\end{eqnarray}
\end{corollary}
\begin{proof}
Using (\ref{2.9})--(\ref{2.12}), we can easily prove that (\ref{2.13})--(\ref{2.15}).
\end{proof}
\begin{theorem}
For $\ell\in \Bbb{N}$, the $q$-recursion formulas of $\mathbf{\Phi}_{1}$, $\mathbf{\Phi}_{2}$ and $\mathbf{\Phi}_{3}$ with the numerator parameters hold true
\begin{eqnarray}
\begin{split}
\mathbf{\Phi}_{1}(q^{\alpha+\ell})=&\mathbf{\Phi}_{1}(q^{\alpha},q^{\beta};q^{\gamma};q,x,y)+\frac{q^{\alpha}x(1-q^{\beta})}{1-q^{\gamma}}\sum_{r=1}^{\ell}q^{r-1}\mathbf{\Phi}_{1}(q^{\alpha+r},q^{\beta+1};q^{\gamma+1};q,x,y)\\
&+\frac{q^{\alpha}y}{1-q^{\gamma}}\sum_{r=1}^{\ell}q^{r-1}\mathbf{\Phi}_{1}(q^{\alpha+r},q^{\beta};q^{\gamma+1};q,qx,y),q^{\gamma}\neq 1,\label{2.16}
\end{split}
\end{eqnarray}\vspace{-12pt}
\begin{eqnarray}
\begin{split}
\mathbf{\Phi}_{1}(q^{\alpha+\ell})=&\mathbf{\Phi}_{1}(q^{\alpha},q^{\beta};q^{\gamma};q,x,y)+\frac{q^{\alpha}y}{1-q^{\gamma}}\sum_{r=1}^{\ell}q^{r-1}\mathbf{\Phi}_{1}(q^{\alpha+r},q^{\beta+1};q^{\gamma+1};q,x,y)\\
&+\frac{q^{\alpha}x(1-q^{\beta})}{1-q^{\gamma}}\sum_{r=1}^{\ell}q^{r-1}\mathbf{\Phi}_{1}(q^{\alpha+r},q^{\beta};q^{\gamma+1};q,x,qy),q^{\gamma}\neq 1,\\
\mathbf{\Phi}_{1}(q^{\beta+\ell})=&\mathbf{\Phi}_{1}(q^{\alpha},q^{\beta};q^{\gamma};q,x,y)+\frac{q^{\beta}x(1-q^{\alpha})}{1-q^{\gamma}}\sum_{r=1}^{\ell}q^{r-1}\mathbf{\Phi}_{1}(q^{\alpha},q^{\beta+\ell};q^{\gamma+1};q,x,y),\label{2.17}
\end{split}
\end{eqnarray}\vspace{-12pt}
\begin{eqnarray}
\begin{split}
\mathbf{\Phi}_{2}(q^{\alpha+\ell})=&\mathbf{\Phi}_{2}(q^{\alpha},q^{\beta};q^{\gamma};q,x,y)+\frac{q^{\alpha}x}{1-q^{\gamma}}\sum_{r=1}^{\ell}q^{r-1}\mathbf{\Phi}_{2}(q^{\alpha+r};q^{\gamma+1};q,x,y),q^{\gamma}\neq 1,\\
\mathbf{\Phi}_{2}(q^{\beta+\ell})=&\mathbf{\Phi}_{2}(q^{\alpha},q^{\beta};q^{\gamma};q,x,y)+\frac{q^{\alpha}y}{1-q^{\gamma}}\sum_{r=1}^{\ell}q^{r-1}\mathbf{\Phi}_{2}(q^{\beta+r};q^{\gamma+1};q,qx,y),q^{\gamma}\neq 1\label{2.18}
\end{split}
\end{eqnarray}
and
\begin{eqnarray}
\begin{split}
\mathbf{\Phi}_{3}(q^{\alpha+\ell})=&\mathbf{\Phi}_{3}(q^{\alpha};q^{\beta};q,x,y)+\frac{q^{\alpha}x}{1-q^{\beta}}\sum_{r=1}^{\ell}q^{r-1}\mathbf{\Phi}_{3}(q^{\alpha+r};q^{\beta+1};q,x,y),q^{\beta}\neq 1.\label{2.19}
\end{split}
\end{eqnarray}
\end{theorem}
\begin{proof}
From (\ref{1.4}) and (\ref{1.7}), we obtain the $q$-contiguous relation:
\begin{eqnarray}
\begin{split}
\mathbf{\Phi}_{1}(q^{\alpha+1})=&\sum_{n,k=0}^{\infty}\bigg{[}1+q^{\alpha}\frac{1-q^{n}}{1-q^{\alpha}}+q^{\alpha+n}\frac{1-q^{k}}{1-q^{\alpha}}\bigg{]}\frac{(q^{\alpha};q)_{n+k}(q^{\beta};q)_{n}x^{n}y^{k}}{(q^{\gamma};q)_{n+k}(q;q)_{n}(q;q)_{k}}\\
=&\mathbf{\Phi}_{1}(q^{\alpha},q^{\beta};q^{\gamma};q,x,y)+\frac{q^{\alpha}x(1-q^{\beta})}{1-q^{\gamma}}\mathbf{\Phi}_{1}(q^{\alpha+1},q^{\beta+1};q^{\gamma+1};q,x,y)\\
&+\frac{q^{\alpha}y}{1-q^{\gamma}}\mathbf{\Phi}_{1}(q^{\alpha+1};q^{\gamma+1};q,qx,y),q^{\gamma}\neq1.\label{2.20}
\end{split}
\end{eqnarray}

Repeating this computation on $\mathbf{\Phi}_{1}$ for $\ell$-times, we obtain the $q$-recursion formula (\ref{2.16}) with parameter $q^{\alpha+\ell}$. Similarly, we get the $q$-recursion formulas (\ref{2.17})--(\ref{2.19}).
\end{proof}
\begin{theorem} The $q$-recursion formulas of $\mathbf{\Phi}_{1}$, $\mathbf{\Phi}_{2}$ and $\mathbf{\Phi}_{3}$ with the lower or denominator parameter hold true
\begin{eqnarray}
\begin{split}
&\mathbf{\Phi}_{1}(q^{\gamma-\ell})=\mathbf{\Phi}_{1}+q^{\gamma}x(1-q^{\alpha})(1-q^{\beta})\sum_{r=1}^{\ell}\frac{q^{r-1}}{(q^{r}-q^{\gamma})(q^{r-1}-q^{\gamma})}\mathbf{\Phi}_{1}(q^{\alpha+1},q^{\beta+1};q^{\gamma+2-r};q,x,y)\\
&+q^{\gamma}y(1-q^{\alpha})\sum_{r=1}^{\ell}\frac{q^{r-1}}{(q^{r}-q^{\gamma})(q^{r-1}-q^{\gamma})}\mathbf{\Phi}_{1}(q^{\alpha+1},q^{\beta};q^{\gamma+2-r};q,q^{r}x,y),q^{\gamma}\neq q^{r},q^{r-1}, r\in \Bbb{N},\label{2.21}
\end{split}
\end{eqnarray}
\begin{eqnarray}
\begin{split}
&\mathbf{\Phi}_{1}(q^{\gamma-\ell})=\mathbf{\Phi}_{1}+q^{\gamma}y(1-q^{\alpha})\sum_{r=1}^{\ell}\frac{q^{r-1}}{(q^{r}-q^{\gamma})(q^{r-1}-q^{\gamma})}\mathbf{\Phi}_{1}(q^{\alpha+1},q^{\beta+1};q^{\gamma+2-r};q,x,y)\\
&+q^{\gamma}x(1-q^{\alpha})(1-q^{\beta})\sum_{r=1}^{\ell}\frac{q^{r-1}}{(q^{r}-q^{\gamma})(q^{r-1}-q^{\gamma})}\mathbf{\Phi}_{1}(q^{\alpha+1},q^{\beta};q^{\gamma+2-r};q,x,q^{r}y),q^{\gamma}\neq q^{r},q^{r-1}, r\in \Bbb{N},\label{2.22}
\end{split}
\end{eqnarray}
\begin{eqnarray}
\begin{split}
\mathbf{\Phi}_{2}(q^{\gamma-\ell})&=\mathbf{\Phi}_{2}+q^{\gamma}x(1-q^{\alpha})\sum_{r=1}^{\ell}\frac{q^{r-1}}{(q^{r}-q^{\gamma})(q^{r-1}-q^{\gamma})}\mathbf{\Phi}_{2}(q^{\alpha+1};q^{\gamma+2-r};q,x,y)\\
&+q^{\gamma}y(1-q^{\beta})\sum_{r=1}^{\ell}\frac{q^{r-1}}{(q^{r}-q^{\gamma})(q^{r-1}-q^{\gamma})}\mathbf{\Phi}_{2}(q^{\beta+1};q^{\gamma+2-r};q,q^{r}x,y),q^{\gamma}\neq q^{r},q^{r-1}, r\in \Bbb{N},\\
\mathbf{\Phi}_{2}(q^{\gamma-\ell})&=\mathbf{\Phi}_{2}+q^{\gamma}y(1-q^{\beta})\sum_{r=1}^{\ell}\frac{q^{r-1}}{(q^{r}-q^{\gamma})(q^{r-1}-q^{\gamma})}\mathbf{\Phi}_{2}(q^{\beta+1};q^{\gamma+2-r};q,x,y)\\
&+q^{\gamma}x(1-q^{\alpha})\sum_{r=1}^{\ell}\frac{q^{r-1}}{(q^{r}-q^{\gamma})(q^{r-1}-q^{\gamma})}\mathbf{\Phi}_{2}(q^{\alpha+1};q^{\gamma+2-r};q,x,q^{r}y),q^{\gamma}\neq q^{r},q^{r-1}, r\in \Bbb{N}\label{2.23}
\end{split}
\end{eqnarray}
and
\begin{eqnarray}
\begin{split}
\mathbf{\Phi}_{3}(q^{\beta-\ell})&=\mathbf{\Phi}_{3}+q^{\beta}x(1-q^{\alpha})\sum_{r=1}^{n}\frac{q^{r-1}}{(q^{r}-q^{\beta})(q^{r-1}-q^{\beta})}\mathbf{\Phi}_{3}(q^{\alpha+1};q^{\beta+2-r};q,x,y)\\
&+q^{\beta}y\sum_{r=1}^{n}\frac{q^{r-1}}{(q^{r}-q^{\beta})(q^{r-1}-q^{\beta})}\mathbf{\Phi}_{3}(-;q^{\beta+2-r};q,q^{r}x,y),q^{\beta}\neq q^{r},q^{r-1}, r\in \Bbb{N},\\
\mathbf{\Phi}_{3}(q^{\beta-\ell})&=\mathbf{\Phi}_{3}+q^{\beta}y\sum_{r=1}^{n}\frac{q^{r-1}}{(q^{r}-q^{\beta})(q^{r-1}-q^{\beta})}\mathbf{\Phi}_{3}(-;q^{\beta+2-r};q,x,y)\\
&+q^{\beta}x(1-q^{\alpha})\sum_{r=1}^{n}\frac{q^{r-1}}{(q^{r}-q^{\beta})(q^{r-1}-q^{\beta})}\mathbf{\Phi}_{3}(q^{\alpha+1};q^{\beta+2-r};q,x,q^{r}y),q^{\beta}\neq q^{r},q^{r-1}, r\in \Bbb{N}.\label{2.24}
\end{split}
\end{eqnarray}
\end{theorem}
\begin{proof}
From (\ref{1.5}) and (\ref{1.7}), we get the $q$-contiguous relation:
\begin{eqnarray}
\begin{split}
\mathbf{\Phi}_{1}(q^{\gamma-1})=&\sum_{n,k=0}^{\infty}\bigg{(}1+q^{\gamma-1}\frac{1-q^{n}}{1-q^{\gamma-1}}+q^{\gamma+n-1}\frac{1-q^{k}}{1-q^{\gamma-1}}\bigg{)}\frac{(q^{\alpha};q)_{n+k}(q^{\beta};q)_{n}x^{n}y^{k}}{(q^{\gamma};q)_{n+k}(q;q)_{k}(q;q)_{n}},\\
=&\mathbf{\Phi}_{1}(q^{\alpha},q^{\beta};q^{\gamma};q,x,y)+\frac{q^{\gamma}x(1-q^{\alpha})(1-q^{\beta})}{(q-q^{\gamma})(1-q^{\gamma})}\mathbf{\Phi}_{1}(q^{\alpha+1},q^{\beta+1};q^{\gamma+1};q,x,y)\\
&+\frac{q^{\gamma}y(1-q^{\alpha})}{(q-q^{\gamma})(1-q^{\gamma})}\mathbf{\Phi}_{1}(q^{\alpha+1};q^{\gamma+1};q,qx,y),q^{\gamma}\neq 1,q^{\gamma}\neq q.\label{2.25}
\end{split}
\end{eqnarray}

Repeating this computation on $\mathbf{\Phi}_{1}$ for $\ell-$times, we obtain the $q$-recursion formula (\ref{2.21}) with parameter $q^{\gamma-\ell}$. Similarly, we obtain (\ref{2.22})--(\ref{2.24}).
\end{proof}

\begin{theorem}
The functions $\mathbf{\Phi}_{1}$, $\mathbf{\Phi}_{2}$ and $\mathbf{\Phi}_{3}$ with respect to parameters satisfy the $q$-difference equations\vspace{-6pt}
\begin{eqnarray}
\begin{split}
\mathbb{D}_{\alpha,q}\mathbf{\Phi}_{1}=&-\frac{1}{1-q^{\alpha}}\bigg{[}\mathbf{\Theta}_{x,q}\mathbf{\Phi}_{1}+\mathbf{\Theta}_{y,q}\mathbf{\Phi}_{1}(qx)\bigg{]},q^{\alpha}\neq1,\label{2.26}
\end{split}
\end{eqnarray}\vspace{-12pt}
\begin{eqnarray}
\begin{split}
\mathbb{D}_{\alpha,q}\mathbf{\Phi}_{1}=&-\frac{1}{1-q^{\alpha}}\bigg{[}\mathbf{\Theta}_{y,q}\mathbf{\Phi}_{1}+\mathbf{\Theta}_{x,q}\mathbf{\Phi}_{1}(qy)\bigg{]},q^{\alpha}\neq1,\\
\mathbb{D}_{\beta,q}\mathbf{\Phi}_{1}=&-\frac{1}{1-q^{\beta}}\mathbf{\Theta}_{x,q}\mathbf{\Phi}_{1},q^{\beta}\neq1,\\
\mathbb{D}_{\gamma,q}\mathbf{\Phi}_{1}=&\frac{1}{1-q^{\gamma}}\bigg{[}\mathbf{\Theta}_{x,q}\mathbf{\Phi}_{1}(q^{\gamma+1})+\mathbf{\Theta}_{y,q}\mathbf{\Phi}_{1}(q^{\gamma+1},qx)\bigg{]},q^{\gamma}\neq1,\\
\mathbb{D}_{\gamma,q}\mathbf{\Phi}_{1}=&\frac{1}{1-q^{\gamma}}\bigg{[}\mathbf{\Theta}_{y,q}\mathbf{\Phi}_{1}(q^{\gamma+1})+\mathbf{\Theta}_{x,q}\mathbf{\Phi}_{1}(q^{\gamma+1},qy)\bigg{]},q^{\gamma}\neq1,\label{2.27}
\end{split}
\end{eqnarray}\vspace{-12pt}
\begin{eqnarray}
\begin{split}
\mathbb{D}_{\alpha,q}\mathbf{\Phi}_{2}=&-\frac{1}{1-q^{\alpha}}\mathbf{\Theta}_{x,q}\mathbf{\Phi}_{2},q^{\alpha}\neq1,\\
\mathbb{D}_{\beta,q}\mathbf{\Phi}_{2}=&-\frac{1}{1-q^{\beta}}\mathbf{\Theta}_{y,q}\mathbf{\Phi}_{2},q^{\beta}\neq1,\\
\mathbb{D}_{\gamma,q}\mathbf{\Phi}_{2}=&\frac{1}{1-q^{\gamma}}\bigg{[}\mathbf{\Theta}_{x,q}\mathbf{\Phi}_{2}(q^{\gamma+1})+\mathbf{\Theta}_{y,q}\mathbf{\Phi}_{2}(q^{\gamma+1},qx)\bigg{]},q^{\gamma}\neq1,\\
\mathbb{D}_{\gamma,q}\mathbf{\Phi}_{2}=&\frac{1}{1-q^{\gamma}}\bigg{[}\mathbf{\Theta}_{y,q}\mathbf{\Phi}_{2}(q^{\gamma+1})+\mathbf{\Theta}_{x,q}\mathbf{\Phi}_{2}(q^{\gamma+1},qy)\bigg{]},q^{\gamma}\neq1\label{2.28}
\end{split}
\end{eqnarray}\vspace{-12pt}
and
\begin{eqnarray}
\begin{split}
\mathbb{D}_{\alpha,q}\mathbf{\Phi}_{3}=&-\frac{1}{1-q^{\alpha}}\mathbf{\Theta}_{x,q}\mathbf{\Phi}_{3},q^{\alpha}\neq1,\\
\mathbb{D}_{\beta,q}\mathbf{\Phi}_{3}=&\frac{1}{1-q^{\beta}}\bigg{[}\mathbf{\Theta}_{x,q}\mathbf{\Phi}_{3}(q^{\beta+1})+\mathbf{\Theta}_{y,q}\mathbf{\Phi}_{3}(q^{\beta+1},qx)\bigg{]},q^{\beta}\neq1,\\
\mathbb{D}_{\beta,q}\mathbf{\Phi}_{3}=&\frac{1}{1-q^{\beta}}\bigg{[}\mathbf{\Theta}_{y,q}\mathbf{\Phi}_{3}(q^{\beta+1})+\mathbf{\Theta}_{x,q}\mathbf{\Phi}_{3}(q^{\beta+1},qy)\bigg{]},q^{\beta}\neq1.\label{2.29}
\end{split}
\end{eqnarray}
\end{theorem}
\begin{proof}
Applying the $q$-difference operator $\mathbb{D}_{\alpha,q}$ (\ref{1.6}) and using (\ref{1.7}), we have
\begin{eqnarray*}
\begin{split}
\mathbb{D}_{\alpha,q}\mathbf{\Phi}_{1}=&\sum_{n,k=0}^{\infty}\frac{(q^{\alpha};q)_{n+k}-(q^{\alpha+1};q)_{n+k}}{(1-q)q^{\alpha}}\frac{(q^{\beta};q)_{n}}{(q^{\gamma};q)_{n+k}}\frac{x^{n}y^{k}}{(q;q)_{k}(q;q)_{n}}\\
=&\sum_{n,k=0}^{\infty}\bigg{(}\frac{q^{\alpha+n+k}-q^{\alpha}}{1-q^{\alpha}}\bigg{)}\frac{(q^{\alpha};q)_{n+k}(q^{\beta};q)_{n}}{(1-q)q^{\alpha}(q^{\gamma};q)_{n+k}}\frac{x^{n}y^{k}}{(q;q)_{k}(q;q)_{n}}\\
=&-\frac{1}{1-q^{\alpha}}\sum_{n,k=0}^{\infty}\bigg{(}\frac{1-q^{n+k}}{1-q}\bigg{)}\frac{(q^{\alpha};q)_{n+k}(q^{\beta};q)_{n}}{(q^{\gamma};q)_{n+k}}\frac{x^{n}y^{k}}{(q;q)_{k}(q;q)_{n}}\\
=&-\frac{1}{1-q^{\alpha}}\sum_{n,k=0}^{\infty}\bigg{(}\frac{1-q^{n}}{1-q}+q^{n}\frac{1-q^{k}}{1-q}\bigg{)}\frac{(q^{\alpha};q)_{n+k}(q^{\beta};q)_{n}}{(q^{\gamma};q)_{n+k}}\frac{x^{n}y^{k}}{(q;q)_{k}(q;q)_{n}}\\
=&-\frac{1}{1-q^{\alpha}}\bigg{[}\mathbf{\Theta}_{x,q}\mathbf{\Phi}_{1}+\mathbf{\Theta}_{y,q}\mathbf{\Phi}_{1}(qx)\bigg{]},q^{\alpha}\neq1.
\end{split}
\end{eqnarray*}
By the same technique in the above, we obtain (\ref{2.27}), (\ref{2.28}) and (\ref{2.29}).
\end{proof}
\begin{theorem}
For the functions $\mathbf{\Phi}_{1}$, $\mathbf{\Phi}_{1}$ and $\mathbf{\Phi}_{3}$, we have the relations
\begin{eqnarray}
\begin{split}
[\mathbf{\Theta}_{x,q}]_{q}\mathbf{\Phi}_{1}=\frac{[\alpha]_{q}[\beta]_{q}}{[\gamma]_{q}}x\mathbf{\Phi}_{1}(q^{\alpha+1},q^{\beta+1};q^{\gamma+1};q,x,y),\label{2.30}
\end{split}
\end{eqnarray}\vspace{-12pt}
\begin{eqnarray}
\begin{split}
[\mathbf{\Theta}_{y,q}]_{q}\mathbf{\Phi}_{1}=\frac{[\alpha]_{q}}{(1-q)[\gamma]_{q}}y\mathbf{\Phi}_{1}(q^{\alpha+1},q^{\beta};q^{\gamma+1};q,x,y),\label{2.31}
\end{split}
\end{eqnarray}\vspace{-12pt}
\begin{eqnarray}
\begin{split}
[\mathbf{\Theta}_{x,q}]_{q}\mathbf{\Phi}_{2}=&\frac{[\alpha]_{q}}{(1-q)[\gamma]_{q}}x\mathbf{\Phi}_{2}(q^{\alpha+1};q^{\gamma+1};q,x,y),\\
[\mathbf{\Theta}_{y,q}]_{q}\mathbf{\Phi}_{2}=&\frac{[\beta]_{q}}{(1-q)[\gamma]_{q}}y\mathbf{\Phi}_{2}(q^{\beta+1};q^{\gamma+1};q,x,y)\label{2.32}
\end{split}
\end{eqnarray}\vspace{-12pt}
and
\begin{eqnarray}
\begin{split}
[\mathbf{\Theta}_{x,q}]_{q}\mathbf{\Phi}_{3}=&\frac{[\alpha]_{q}}{(1-q)[\beta]_{q}}x\mathbf{\Phi}_{3}(q^{\alpha+1};q^{\beta+1};q,x,y),\\
[\mathbf{\Theta}_{y,q}]_{q}\mathbf{\Phi}_{3}=&\frac{1}{(1-q)^{2}[\beta]_{q}}y\mathbf{\Phi}_{3}(q^{\beta+1};q,x,y).\label{2.33}
\end{split}
\end{eqnarray}
\end{theorem}
\begin{proof}
Applying the operator $\mathbf{\Theta}_{x,q}$ to the both sides of (\ref{1.7}) with respect to $x$, we have
\begin{eqnarray*}
\begin{split}
[\mathbf{\Theta}_{x,q}]_{q}\mathbf{\Phi}_{1}=&\sum_{n,k=0}^{\infty}\bigg{(}\frac{1-q^{n}}{1-q}\bigg{)}\frac{(q^{\alpha};q)_{n+k}(q^{\beta};q)_{n}}{(q^{\gamma};q)_{n+k}}\frac{x^{n}y^{k}}{(q;q)_{k}(q;q)_{n}}\\
=&\frac{1}{1-q}\sum_{n=1,k=0}^{\infty}\frac{(q^{\alpha};q)_{n+k}(q^{\beta};q)_{n}x^{n}y^{k}}{(q^{\gamma};q)_{n+k}(q;q)_{n-1}(q;q)_{k}}\\
=&\frac{(1-q^{\alpha})(1-q^{\beta})}{(1-q)(1-q^{\gamma})}\sum_{n,k=0}^{\infty}\frac{(q^{\alpha+1};q)_{n+k}(q^{\beta+1};q)_{n}x^{n+1}y^{k}}{(q^{\gamma+1};q)_{n+k}(q;q)_{k}(q;q)_{n}}\\
=&\frac{(1-q^{\alpha})(1-q^{\beta})}{(1-q^{\gamma})(1-q)}x\mathbf{\Phi}_{1}(q^{\alpha+1},q^{\beta+1};q^{\gamma+1};q,x,y),q^{\gamma}\neq1.
\end{split}
\end{eqnarray*}

By the same way, the proof of Equations (\ref{2.31})--(\ref{2.33}) are similar lines to the proof of Equation (\ref{2.30}).
\end{proof}
\begin{theorem}
The functions $\mathbf{\Phi}_{1}$, $\mathbf{\Phi}_{2}$ and $\mathbf{\Phi}_{3}$ satisfies the $q$-differential relations
\begin{eqnarray}
\begin{split}
[\mathbf{\Theta}_{x,q}+\mathbf{\Theta}_{y,q}+\alpha]_{q}\mathbf{\Phi}_{1}=[\alpha]_{q}\mathbf{\Phi}_{1}(q^{\alpha+1}),\label{2.34}
\end{split}
\end{eqnarray}\vspace{-12pt}
\begin{eqnarray}
\begin{split}
&[\mathbf{\Theta}_{x,q}+\beta]_{q}\mathbf{\Phi}_{1}=[\beta]_{q}\mathbf{\Phi}_{1}(q^{\beta+1}),\\
&[\mathbf{\Theta}_{x,q}+\mathbf{\Theta}_{y,q}+\gamma-1]_{q}\mathbf{\Phi}_{1}=[\gamma-1]_{q}\mathbf{\Phi}_{1}(q^{\gamma-1}),\label{2.35}
\end{split}
\end{eqnarray}\vspace{-12pt}
\begin{eqnarray}
\begin{split}
&[\mathbf{\Theta}_{x,q}+\alpha]_{q}\mathbf{\Phi}_{2}=[\alpha]_{q}\mathbf{\Phi}_{2}(q^{\alpha+1}),\\
&[\mathbf{\Theta}_{y,q}+\beta]_{q}\mathbf{\Phi}_{2}=[\beta]_{q}\mathbf{\Phi}_{2}(q^{\beta+1}),\\
&[\mathbf{\Theta}_{x,q}+\mathbf{\Theta}_{y,q}+\gamma-1]_{q}\mathbf{\Phi}_{2}=[\gamma-1]_{q}\mathbf{\Phi}_{2}(q^{\gamma-1})\label{2.36}
\end{split}
\end{eqnarray}
and
\begin{eqnarray}
\begin{split}
&[\mathbf{\Theta}_{x,q}+\alpha]_{q}\mathbf{\Phi}_{3}=[\alpha]_{q}\mathbf{\Phi}_{3}(q^{\alpha+1}),\\
&[\mathbf{\Theta}_{x,q}+\mathbf{\Theta}_{y,q}+\beta-1]_{q}\mathbf{\Phi}_{3}=[\beta-1]_{q}\mathbf{\Phi}_{3}(q^{\beta-1}).\label{2.37}
\end{split}
\end{eqnarray}
\end{theorem}
\begin{proof}
Using the relation (\ref{1.3}) and applying the $q$-derivatives operators (\ref{1.6}) and (\ref{1.8}), we get
\begin{eqnarray*}
\begin{split}
[\mathbf{\Theta}_{x,q}+\mathbf{\Theta}_{y,q}+\alpha]_{q}\mathbf{\Phi}_{1}=&\sum_{n,k=0}^{\infty}\frac{[\alpha+n+k]_{q}(q^{\alpha};q)_{n+k}(q^{\beta};q)_{n}}{(q^{\gamma};q)_{n+k}(q;q)_{k}(q;q)_{n}}x^{n}y^{k}\\
=&[\alpha]_{q}\sum_{n,k=0}^{\infty}\frac{(q^{\alpha+1};q)_{n+k}(q^{\beta};q)_{n}}{(q^{\gamma};q)_{n+k}(q;q)_{k}(q;q)_{n}}x^{n}y^{k}=[\alpha]_{q}\mathbf{\Phi}_{1}(q^{\alpha+1}).
\end{split}
\end{eqnarray*}
Using the same procedure in the proof of Equation (\ref{2.34}) leads to the results (\ref{2.35})--(\ref{2.37}).
\end{proof}
\begin{theorem}
The functions $\mathbf{\Phi}_{1}$, $\mathbf{\Phi}_{2}$, and $\mathbf{\Phi}_{3}$ satisfy the $q$-recurrence relations
\begin{eqnarray}
\begin{split}
(1-q^{\alpha})\mathbf{\Phi}_{1}(q^{\alpha+1})+q^{\alpha}\mathbf{\Phi}_{1}(xq,yq)-\mathbf{\Phi}_{1}=0,\label{2.38}
\end{split}
\end{eqnarray}\vspace{-12pt}
\begin{eqnarray}
\begin{split}
&(1-q^{\beta})\mathbf{\Phi}_{1}(q^{\beta+1})+q^{\beta}\mathbf{\Phi}_{1}(xq)-\mathbf{\Phi}_{1}=0,\\
&(1-q^{\gamma-1})\mathbf{\Phi}_{1}(q^{\gamma-1})+q^{\gamma-1}\mathbf{\Phi}_{1}(xq,yq)-\mathbf{\Phi}_{1}=0,\label{2.39}
\end{split}
\end{eqnarray}
\begin{eqnarray}
\begin{split}
&(1-q^{\alpha})\mathbf{\Phi}_{2}(q^{\alpha+1})+q^{\alpha}\mathbf{\Phi}_{2}(xq)-\mathbf{\Phi}_{2}=0,\\
&(1-q^{\beta})\mathbf{\Phi}_{2}(q^{\beta+1})+q^{\beta}\mathbf{\Phi}_{2}(yq)-\mathbf{\Phi}_{2}=0,\\
&(1-q^{\gamma-1})\mathbf{\Phi}_{2}(q^{\gamma-1})+q^{\gamma-1}\mathbf{\Phi}_{2}(xq,yq)-\mathbf{\Phi}_{2}=0\label{2.40}
\end{split}
\end{eqnarray}\vspace{-12pt}
and
\begin{eqnarray}
\begin{split}
&(1-q^{\alpha})\mathbf{\Phi}_{3}(q^{\alpha+1})+q^{\alpha}\mathbf{\Phi}_{3}(xq)-\mathbf{\Phi}_{3}=0,\\
&(1-q^{\beta-1})\mathbf{\Phi}_{3}(q^{\beta-1})+q^{\beta-1}\mathbf{\Phi}_{3}(xq,yq)-\mathbf{\Phi}_{3}=0.\label{2.41}
\end{split}
\end{eqnarray}
\end{theorem}
\begin{proof}
Applying an operator, we have
\begin{eqnarray*}
\begin{split}
[\mathbf{\Theta}_{x,q}+\mathbf{\Theta}_{y,q}+\alpha]_{q}\mathbf{\Phi}_{1}=\frac{1-q^{\mathbf{\Theta}_{x,q}+\mathbf{\Theta}_{y,q}+\alpha}}{1-q}\mathbf{\Phi}_{1}\\
=\frac{\mathbf{\Phi}_{1}-q^{\alpha}\mathbf{\Phi}_{1}(xq,yq)}{1-q}=[\alpha]_{q}\mathbf{\Phi}_{1}(q^{\alpha+1})
\end{split}
\end{eqnarray*}
and using (\ref{2.34}), we obtain (\ref{2.38}). The proofs of (\ref{2.39})--(\ref{2.41}) follow in the same way.
\end{proof}
\begin{theorem}
The functions $\mathbf{\Phi}_{1}$, $\mathbf{\Phi}_{2}$ and $\mathbf{\Phi}_{3}$ satisfy the $q$-partial differential equations
\begin{eqnarray}
\begin{split}
\bigg{(}[\mathbf{\Theta}_{x,q}]_{q}[\mathbf{\Theta}_{x,q}+\mathbf{\Theta}_{y,q}+\gamma-1]_{q}-x[\mathbf{\Theta}_{x,q}+\mathbf{\Theta}_{y,q}+\alpha]_{q}[\mathbf{\Theta}_{x,q}+\beta]_{q}\bigg{)}\mathbf{\Phi}_{1}=0,\label{2.42}
\end{split}
\end{eqnarray}\vspace{-12pt}
\begin{eqnarray}
\begin{split}
\bigg{(}[\mathbf{\Theta}_{y,q}]_{q}[\mathbf{\Theta}_{x,q}+\mathbf{\Theta}_{y,q}+\gamma-1]_{q}-\frac{y}{1-q}[\mathbf{\Theta}_{x,q}+\mathbf{\Theta}_{y,q}+\alpha]_{q}\bigg{)}\mathbf{\Phi}_{1}=0,\label{2.43}
\end{split}
\end{eqnarray}
\begin{eqnarray}
\begin{split}
\bigg{(}[\mathbf{\Theta}_{x,q}]_{q}[\mathbf{\Theta}_{x,q}+\mathbf{\Theta}_{y,q}+\gamma-1]_{q}-\frac{x}{1-q}[\mathbf{\Theta}_{x,q}+\alpha]_{q}\bigg{)}\mathbf{\Phi}_{2}=0,\\
\bigg{(}[\mathbf{\Theta}_{y,q}]_{q}[\mathbf{\Theta}_{x,q}+\mathbf{\Theta}_{y,q}+\gamma-1]_{q}-\frac{y}{1-q}[\mathbf{\Theta}_{y,q}+\beta]_{q}\bigg{)}\mathbf{\Phi}_{2}=0\label{2.44}
\end{split}
\end{eqnarray}
and
\begin{eqnarray}
\begin{split}
\bigg{(}[\mathbf{\Theta}_{x,q}]_{q}[\mathbf{\Theta}_{x,q}+\mathbf{\Theta}_{y,q}+\beta-1]_{q}-\frac{x}{1-q}[\mathbf{\Theta}_{x,q}+\alpha]_{q}\bigg{)}\mathbf{\Phi}_{3}=0,\\
\bigg{(}[\mathbf{\Theta}_{y,q}]_{q}[\mathbf{\Theta}_{x,q}+\mathbf{\Theta}_{y,q}+\beta-1]_{q}-\frac{y}{(1-q)^{2}}\bigg{)}\mathbf{\Phi}_{3}=0.\label{2.45}
\end{split}
\end{eqnarray}
\end{theorem}
\begin{proof}
Applying the $q$-derivatives operators (\ref{1.6}) and using (\ref{2.30}), (\ref{2.34}) and (\ref{2.35}), we get
\begin{eqnarray*}
\begin{split}
&[\mathbf{\Theta}_{x,q}]_{q}[\mathbf{\Theta}_{x,q}+\mathbf{\Theta}_{y,q}+\gamma-1]_{q}\mathbf{\Phi}_{1}=\sum_{n,k=0}^{\infty}\frac{[n]_{q}[n+k+\gamma-1]_{q}(q^{\alpha};q)_{n+k}(q^{\beta};q)_{n}}{(q^{\gamma};q)_{n+k}}\frac{x^{n}y^{k}}{(q;q)_{k}(q;q)_{n}}\\
&=\sum_{n=1,k=0}^{\infty}\frac{(q^{\alpha};q)_{n+k}(q^{\beta};q)_{n}}{(q^{\gamma};q)_{n+k-1}}\frac{x^{n}y^{k}}{(q;q)_{n-1}(q;q)_{k}}\\
&=\sum_{n,k=0}^{\infty}\frac{[\alpha+n+k]_{q}[\beta+n]_{q}(q^{\alpha};q)_{n+k}(q^{\beta};q)_{n}}{(q^{\gamma};q)_{n+k}}\frac{x^{n+1}y^{k}}{(q;q)_{k}(q;q)_{n}}\\
&=x[\mathbf{\Theta}_{x,q}+\mathbf{\Theta}_{y,q}+\alpha]_{q}[\mathbf{\Theta}_{x,q}+\beta]_{q}\mathbf{\Phi}_{1}.
\end{split}
\end{eqnarray*}
Same as, we obtain the $q$-partial differential Equations (\ref{2.43})--(\ref{2.45}).
\end{proof}
\begin{corollary} The $q$-partial differential equations hold true for the functions $\mathbf{\Phi}_{1}$, $\mathbf{\Phi}_{2}$ and $\mathbf{\Phi}_{3}$
\begin{eqnarray}
\begin{split}
&\bigg{(}q^{\gamma-1}[\mathbf{\Theta}_{x,q}]_{q}[\mathbf{\Theta}_{x,q}+\mathbf{\Theta}_{y,q}]_{q}-q^{\gamma-1}[\mathbf{\Theta}_{x,q}]_{q}+[\gamma]_{q}[\mathbf{\Theta}_{x,q}]_{q}-xq^{\alpha+\beta}[\mathbf{\Theta}_{x,q}+\mathbf{\Theta}_{y,q}]_{q}[\mathbf{\Theta}_{x,q}]_{q}\\
&-xq^{\alpha}[\beta]_{q}[\mathbf{\Theta}_{x,q}+\mathbf{\Theta}_{y,q}]_{q}-xq^{\beta}[\alpha]_{q}[\mathbf{\Theta}_{x,q}]_{q}-x[\alpha]_{q}[\beta]_{q}\bigg{)}\mathbf{\Phi}_{1}=0,\label{2.46}
\end{split}
\end{eqnarray}
\begin{eqnarray}
\begin{split}
&\bigg{(}q^{\gamma-1}[\mathbf{\Theta}_{y,q}]_{q}[\mathbf{\Theta}_{x,q}+\mathbf{\Theta}_{y,q}]_{q}-\frac{y}{1-q}q^{\alpha}[\mathbf{\Theta}_{x,q}+\mathbf{\Theta}_{y,q}]_{q}+\bigg{(}[\gamma]_{q}-q^{\gamma-1}\bigg{)}[\mathbf{\Theta}_{y,q}]_{q}-\frac{y}{1-q}[\alpha]_{q}\bigg{)}\mathbf{\Phi}_{1}=0,\label{2.47}
\end{split}
\end{eqnarray}
\begin{eqnarray}
\begin{split}
&\bigg{(}q^{\gamma-1}[\mathbf{\Theta}_{x,q}]_{q}[\mathbf{\Theta}_{x,q}+\mathbf{\Theta}_{y,q}]_{q}+\bigg{(}[\gamma]_{q}-q^{\gamma-1}-\frac{x}{1-q}q^{\alpha}\bigg{)}[\mathbf{\Theta}_{x,q}]_{q}-\frac{x}{1-q}[\alpha]_{q}\bigg{)}\mathbf{\Phi}_{2}=0,\\
&\bigg{(}q^{\gamma-1}[\mathbf{\Theta}_{y,q}]_{q}[\mathbf{\Theta}_{x,q}+\mathbf{\Theta}_{y,q}]_{q}+\bigg{(}[\gamma]_{q}-q^{\gamma-1}-\frac{y}{1-q}q^{\beta}\bigg{)}[\mathbf{\Theta}_{y,q}]_{q}-\frac{y}{1-q}[\beta]_{q}\bigg{)}\mathbf{\Phi}_{2}=0\label{2.48}
\end{split}
\end{eqnarray}
and
\begin{eqnarray}
\begin{split}
&\bigg{(}q^{\beta-1}[\mathbf{\Theta}_{x,q}]_{q}[\mathbf{\Theta}_{x,q}+\mathbf{\Theta}_{y,q}]_{q}+\bigg{(}[\beta]_{q}-q^{\beta-1}[1]_{q}-\frac{x}{1-q}q^{\alpha}\bigg{)}[\mathbf{\Theta}_{x,q}]_{q}-\frac{x}{1-q}[\alpha]_{q}\bigg{)}\mathbf{\Phi}_{3}=0,\\
&\bigg{(}q^{\beta-1}[\mathbf{\Theta}_{y,q}]_{q}[\mathbf{\Theta}_{x,q}+\mathbf{\Theta}_{y,q}]_{q}+\bigg{(}[\beta]_{q}-q^{\beta-1}[1]_{q}\bigg{)}[\mathbf{\Theta}_{y,q}]_{q}-\frac{y}{(1-q)^{2}}\bigg{)}\mathbf{\Phi}_{3}=0.\label{2.49}
\end{split}
\end{eqnarray}
\end{corollary}
\begin{proof}
Applying the identity
\begin{eqnarray*}
\begin{split}
[n-k]_{q}&=q^{-k}\bigg{(}[n]_{q}-[k]_{q}\bigg{)},\\
[n+k]_{q}&=[n]_{q}+q^{n}[k]_{q}=[k]_{q}+q^{k}[n]_{q}
\end{split}
\end{eqnarray*}
and writing the $q$-partial differential Equations (\ref{2.35})--(\ref{2.38}), we obtain (\ref{2.46})--(\ref{2.49}).
\end{proof}
\begin{theorem} The summation formulas for the functions $\mathbf{\Phi}_{1}$, $\mathbf{\Phi}_{2}$ and $\mathbf{\Phi}_{3}$ hold true:
\begin{eqnarray}
\begin{split}
\mathbf{\Phi}_{1}(q^{\alpha},q^{\beta};q^{\gamma};q,q,x,y)=\sum_{n=0}^{\infty}\frac{(q^{\alpha};q)_{n}(q^{\beta};q)_{n}}{(q^{\gamma};q)_{n}(q;q)_{n}}x^{n}
\;_{2}\mathbf{\Phi}_{1}(q^{\alpha+n},0;q^{\gamma+n};q,y),\label{2.50}
\end{split}
\end{eqnarray}\vspace{-12pt}
\begin{eqnarray}
\begin{split}
&\mathbf{\Phi}_{1}(q^{\alpha},q^{\beta};q^{\gamma};q,q,x,y)=\sum_{k=0}^{\infty}\frac{(q^{\alpha};q)_{k}}{(q^{\gamma};q)_{k}(q;q)_{k}}y^{k}\;_{2}\mathbf{\Phi}_{1}(q^{\alpha+k},q^{\beta};q^{\gamma+k};q,x),\label{2.51}
\end{split}
\end{eqnarray}\vspace{-12pt}
\begin{eqnarray}
\begin{split}
&\mathbf{\Phi}_{2}(q^{\alpha},q^{\beta};q^{\gamma};q,q,x,y)=\sum_{n=0}^{\infty}\frac{(q^{\alpha};q)_{n}}{(q^{\gamma};q)_{n}(q;q)_{n}}x^{n}
\;_{2}\mathbf{\Phi}_{1}(q^{\beta},0;q^{\gamma+n};q,y),\\
&\mathbf{\Phi}_{2}(q^{\alpha},q^{\beta};q^{\gamma};q,q,x,y)=\sum_{k=0}^{\infty}\frac{(q^{\beta};q)_{k}}{(q^{\gamma};q)_{k}(q;q)_{k}}y^{k}\;_{2}\mathbf{\Phi}_{1}(q^{\alpha},0;q^{\gamma+k};q,x)\label{2.52}
\end{split}
\end{eqnarray}\vspace{-12pt}
and
\begin{eqnarray}
\begin{split}
&\mathbf{\Phi}_{3}(q^{\alpha};q^{\beta};q,q,x,y)=\sum_{n=0}^{\infty}\frac{(q^{\alpha};q)_{n}}{(q^{\beta};q)_{n}(q;q)_{n}}x^{n}
\;_{2}\mathbf{\Phi}_{1}(0,0;q^{\beta+n};q,y),\\
&\mathbf{\Phi}_{3}(q^{\alpha};q^{\beta};q,q,x,y)=\sum_{k=0}^{\infty}\frac{1}{(q^{\beta};q)_{k}(q;q)_{k}}y^{k}\;_{2}\mathbf{\Phi}_{1}(q^{\alpha},0;q^{\beta+k};q,x).\label{2.53}
\end{split}
\end{eqnarray}
\end{theorem}
\begin{proof}
By using the relation
\begin{eqnarray*}
\begin{split}
(q^{\alpha};q)_{n+k}=(q^{\alpha};q)_{n}(q^{\alpha+n};q)_{k},
\end{split}
\end{eqnarray*}
we have
\begin{eqnarray*}
\begin{split}
&\mathbf{\Phi}_{1}(q^{\alpha},q^{\beta};q^{\gamma};q,q,x,y)=\sum_{n,k=0}^{\infty}\frac{(q^{\alpha};q)_{n}(q^{a+n};q)_{k}(q^{\beta};q)_{n}}{(q^{\gamma};q)_{n}(q^{\gamma+n};q)_{k}(q;q)_{k}(q;q)_{n}}x^{n}y^{k}\\
&=\sum_{n=0}^{\infty}\frac{(q^{\alpha};q)_{n}(q^{\beta};q)_{n}}{(q^{\gamma};q)_{n}(q;q)_{n}}x^{n}\sum_{k=0}^{\infty}\frac{(q^{\alpha+n};q)_{k}}{(q^{\gamma+n};q)_{k}(q;q)_{k}}y^{k}\\
&=\sum_{n=0}^{\infty}\frac{(q^{\alpha};q)_{n}(q^{\beta};q)_{n}}{(q^{\gamma};q)_{n}(q;q)_{n}}x^{n}\;_{2}\mathbf{\Phi}_{1}(q^{\alpha+n},0;q^{\gamma+n};q,y).
\end{split}
\end{eqnarray*}

This technique used in (\ref{2.50}) yields new formulas for the basic Humbert hypergeometric $\mathbf{\Phi}_{1}$, $\mathbf{\Phi}_{2}$ and $\mathbf{\Phi}_{3}$ series.
\end{proof}
Some of the interesting special cases of Equations (\ref{2.50})--(\ref{2.53}) are: 
\begin{enumerate}
  \item Setting $y=0$ in (\ref{2.50}), we get
\begin{eqnarray*}
\begin{split}
&\mathbf{\Phi}_{1}(q^{\alpha},q^{\beta};q^{\gamma};q,q,x,0)=\sum_{n=0}^{\infty}\frac{(q^{\alpha};q)_{n}(q^{\beta};q)_{n}}{(q^{\gamma};q)_{n}(q;q)_{n}}x^{n}\;_{2}\mathbf{\Phi}_{1}(q^{\alpha+n},0;q^{\gamma+n};q,0)\\
&=\sum_{n=0}^{\infty}\frac{(q^{\alpha};q)_{n}(q^{\beta};q)_{n}}{(q^{\gamma};q)_{n}(q;q)_{n}}x^{n}=\;_{2}\mathbf{\Phi}_{1}(q^{\alpha},q^{\beta};q^{\gamma};q,x).
\end{split}
\end{eqnarray*}
  \item Taking $x=0$ in (\ref{2.51}), we get
\begin{eqnarray*}
\begin{split}
&\mathbf{\Phi}_{1}(q^{\alpha},q^{\beta};q^{\gamma};q,q,0,y)=\sum_{k=0}^{\infty}\frac{(q^{\alpha};q)_{k}}{(q^{\gamma};q)_{k}(q;q)_{k}}y^{k}\;_{2}\mathbf{\Phi}_{1}(q^{\alpha+k},q^{\beta};q^{\gamma+k};q,0)\\
&=\sum_{k=0}^{\infty}\frac{(q^{\alpha};q)_{k}}{(q^{\gamma};q)_{k}(q;q)_{k}}y^{k}=\;_{2}\mathbf{\Phi}_{1}(q^{\alpha},0;q^{\gamma};q,y).
\end{split}
\end{eqnarray*}
  \item Letting $y=0$ and $x=0$ in (\ref{2.52}), we have
\begin{eqnarray*}
\begin{split}
&\mathbf{\Phi}_{2}(q^{\alpha},q^{\beta};q^{\gamma};q,x,0)=\;_{2}\mathbf{\Phi}_{1}(q^{\alpha},0;q^{\gamma};q,x),\\
&\mathbf{\Phi}_{2}(q^{\alpha},q^{\beta};q^{\gamma};q,0,y)=\;_{2}\mathbf{\Phi}_{1}(q^{\beta},0;q^{\gamma};q,y).
\end{split}
\end{eqnarray*}
  \item Putting $y=0$ and $x=0$ in (\ref{2.53}), we obtain
\begin{eqnarray*}
\begin{split}
&\mathbf{\Phi}_{3}(q^{\alpha};q^{\beta};q,x,0)=\;_{2}\mathbf{\Phi}_{1}(q^{\alpha},0;q^{\beta};q,x),\\
&\mathbf{\Phi}_{3}(q^{\alpha};q^{\beta};q,0,y)=\;_{2}\mathbf{\Phi}_{1}(0,0;q^{\beta};q,y).
\end{split}
\end{eqnarray*}
\end{enumerate}

\begin{theorem}
The transformation formula for $\mathbf{\Phi}_{1}$ is true 
\begin{eqnarray}
\begin{split}
&\mathbf{\Phi}_{1}(q^{\alpha},q^{\beta};q^{\gamma};q,x,y)=\frac{(q^{\alpha},xq^{\beta};q)_{\infty}}{(q^{\gamma},x,y;q)_{\infty}}\;_{3}\mathbf{\Phi}_{2}(q^{\gamma-\alpha},x,y;q^{\beta}x,0;q,q^{\alpha}).\label{2.54}
\end{split}
\end{eqnarray}
\end{theorem}
\begin{proof}
Using the identities\vspace{-6pt}
\begin{eqnarray*}
\begin{split}
(q^{\alpha};q)_{n+k}=\frac{(q^{\alpha};q)_{\infty}}{(q^{\alpha+n+k};q)_{\infty}}\\
(q^{\gamma};q)_{n+k}=\frac{(q^{\gamma};q)_{\infty}}{(q^{\gamma+n+k};q)_{\infty}},
\end{split}
\end{eqnarray*}\vspace{-12pt}
\begin{eqnarray*}
\begin{split}
\frac{(q^{\gamma+n+k};q)_{\infty}}{(q^{\alpha+n+k};q)_{\infty}}=\sum_{s=0}^{\infty}\frac{(q^{\gamma-\alpha};q)_{s}}{(q;q)_{s}}\bigg{(}q^{\alpha+n+k}\bigg{)}^{s},
\end{split}
\end{eqnarray*}
\begin{eqnarray*}
\begin{split}
&\:_{1}\mathbf{\Phi}_{0}(q^{\alpha};-;q,x)=\sum_{n=0}^{\infty}\frac{(q^{\alpha};q)_{n}}{(q;q)_{n}}x^{n}=\frac{(q^{\alpha}x;q)_{\infty}}{(x;q)_{\infty}};|x|<1
\end{split}
\end{eqnarray*}\vspace{-12pt}
and
\begin{eqnarray*}
\begin{split}
&\:_{0}\mathbf{\Phi}_{0}(-;-;q,y)=\sum_{k=0}^{\infty}\frac{1}{(q;q)_{k}}y^{k}=\frac{1}{(y;q)_{\infty}},
\end{split}
\end{eqnarray*}
substitution of this gives 
\begin{eqnarray*}
\begin{split}
&\mathbf{\Phi}_{1}(q^{\alpha},q^{\beta};q^{\gamma};q,x,y)=\frac{(q^{\alpha};q)_{\infty}}{(q^{\gamma};q)_{\infty}}\sum_{n,k=0}^{\infty}\frac{(q^{\gamma+n+k};q)_{\infty}(q^{\beta};q)_{n}}{(q^{\alpha+n+k};q)_{\infty}(q;q)_{n}(q;q)_{k}}x^{n}y^{k}\\
&=\frac{(q^{\alpha};q)_{\infty}}{(q^{\gamma};q)_{\infty}}\sum_{n,k,s=0}^{\infty}\frac{(q^{\gamma-\alpha};q)_{s}(q^{\beta};q)_{n}}{(q;q)_{s}(q;q)_{n}(q;q)_{k}}\bigg{(}q^{\alpha+n+k}\bigg{)}^{s}x^{n}y^{k}\\
&=\frac{(q^{\alpha};q)_{\infty}}{(q^{\gamma};q)_{\infty}}\sum_{s=0}^{\infty}\frac{(q^{\gamma-\alpha};q)_{s}}{(q;q)_{s}}q^{s\alpha}\sum_{n=0}^{\infty}\frac{(q^{\beta};q)_{n}}{(q;q)_{n}}q^{ns}x^{n}
\sum_{k=0}^{\infty}\frac{1}{(q;q)_{k}}q^{ks}y^{k}\\
&=\frac{(q^{\alpha};q)_{\infty}}{(q^{\gamma};q)_{\infty}}\sum_{s=0}^{\infty}\frac{(q^{\gamma-\alpha};q)_{s}}{(q;q)_{s}}q^{s\alpha}\frac{(q^{\beta+s}x;q)_{\infty}}{(q^{s}x;q)_{\infty}}\frac{1}{(q^{s}y;q)_{\infty}}\\
&=\frac{(q^{\alpha};q)_{\infty}(q^{\beta}x;q)_{\infty}}{(q^{\gamma};q)_{\infty}(x;q)_{\infty}(y;q)_{\infty}}\sum_{s=0}^{\infty}\frac{(q^{\gamma-\alpha};q)_{s}(x;q)_{s}(y;q)_{s}}{(q;q)_{s}(q^{\beta}x;q)_{s}}q^{s\alpha}\\
&=\frac{(q^{\alpha};q)_{\infty}(q^{\beta}x;q)_{\infty}}{(q^{\gamma};q)_{\infty}(x;q)_{\infty}(y;q)_{\infty}}\;_{3}\mathbf{\Phi}_{2}(q^{\gamma-\alpha},x,y;q^{\beta}x,0;q,q^{\alpha}).
\end{split}
\end{eqnarray*}
\end{proof}
Some of the interesting special cases of Equation (\ref{2.54}) are: 
\begin{enumerate}
 \item Putting $y=xq^{\beta}$ in (\ref{2.54}), we get
\begin{eqnarray*}
\begin{split}
&\mathbf{\Phi}_{1}(q^{\alpha},q^{\beta};q^{\gamma};q,x,q^{\beta}x)=\frac{(q^{\alpha},q^{\beta}x;q)_{\infty}}{(q^{\gamma},x,q^{\beta}x;q)_{\infty}}\;_{3}\mathbf{\Phi}_{2}(q^{\gamma-\alpha},x,q^{\beta}x;q^{\beta}x,0;q,q^{\alpha})\\
&=\frac{(q^{\alpha};q)_{\infty}}{(q^{\gamma},x;q)_{\infty}}\;_{2}\mathbf{\Phi}_{1}(q^{\gamma-\alpha},x;0;q,q^{\alpha})=\;_{2}\mathbf{\Phi}_{1}(q^{\alpha},0;q^{\gamma};q,x).
\end{split}
\end{eqnarray*}
  \item Letting $x=q^{\gamma-\alpha-\beta}$ in (\ref{2.54}), we have
\begin{eqnarray*}
\begin{split}
&\mathbf{\Phi}_{1}(q^{\alpha},q^{\beta};q^{\gamma};q,q^{\gamma-\alpha-\beta},y)=\frac{(q^{\alpha},q^{\gamma-\alpha-\beta}q^{\beta};q)_{\infty}}{(q^{\gamma},q^{\gamma-\alpha-\beta},y;q)_{\infty}}\;_{3}\mathbf{\Phi}_{2}(q^{\gamma-\alpha},q^{\gamma-\alpha-\beta},y;q^{\beta}q^{\gamma-\alpha-\beta},0;q,q^{\alpha})\\
&=\frac{(q^{\alpha},q^{\gamma-\alpha};q)_{\infty}}{(q^{\gamma},q^{\gamma-\alpha-\beta},y;q)_{\infty}}\;_{3}\mathbf{\Phi}_{2}(q^{\gamma-\alpha},q^{\gamma-\alpha-\beta},y;q^{\gamma-\alpha},0;q,q^{\alpha})\\
&=\frac{(q^{\alpha},q^{\gamma-\alpha};q)_{\infty}}{(q^{\gamma},q^{\gamma-\alpha-\beta},y;q)_{\infty}}\;_{2}\mathbf{\Phi}_{1}(q^{\gamma-\alpha-\beta},y;0;q,q^{\alpha}).
\end{split}
\end{eqnarray*}
\end{enumerate}
\begin{theorem}
For $0<\Re(\alpha)<\Re(\gamma)$, the $q$-integral representations for $\mathbf{\Phi}_{1}$ are true
\begin{eqnarray}
\begin{split}
&\mathbf{\Phi}_{1}(q^{\alpha},q^{\beta};q^{\gamma};q,x,y)=\frac{\Gamma_{q}(\gamma)}{\Gamma_{q}(\alpha)\Gamma_{q}(\gamma-\alpha)}\int_{0}^{1}t^{\alpha-1}\frac{(qt,xtq^{\beta};q)_{\infty}}{(xt,yt,tq^{\gamma-\alpha};q)_{\infty}}d_{q}t\label{2.55}
\end{split}
\end{eqnarray}
where $\Gamma_{q}(\alpha)$ is $q$-gamma function
\begin{eqnarray*}
\begin{split}
\Gamma_{q}(\alpha)=\frac{(q;q)_{\alpha-1}}{(1-q)^{\alpha-1}}=\frac{(q;q)_{\infty}}{(q^{\alpha};q)_{\infty}(1-q)^{\alpha-1}};\alpha\neq 0,-1,-2,\ldots
\end{split}
\end{eqnarray*}
and
\begin{eqnarray}
\begin{split}
&\mathbf{\Phi}_{1}(q^{\alpha},q^{\beta};q^{\gamma};q,x,y)=\frac{\Gamma_{q}(\gamma)}{\Gamma_{q}(\alpha)\Gamma_{q}(\gamma-\alpha)}\int_{0}^{1}t^{\alpha-1}(1-qt)_{\gamma-\alpha-1}(1-qxt)_{-\beta}e_{q}\bigg{(}\frac{ty}{1-q}\bigg{)}d_{q}t,\label{2.56}
\end{split}
\end{eqnarray}
where
\vspace{-12pt}\begin{eqnarray*}
\begin{split}
&(1-qt)_{-\alpha}=1+\frac{[\alpha]_{q}}{[1]_{q}!}t+\frac{[\alpha]_{q}[\alpha+1]_{q}}{[2]_{q}!}t^{2}+\ldots.
\end{split}
\end{eqnarray*}
\end{theorem}
\begin{proof}
Using
\begin{eqnarray}
\begin{split}
\frac{(q^{\alpha};q)_{n+k}}{(q^{\gamma};q)_{n+k}}=\frac{\Gamma_{q}(\gamma)}{\Gamma_{q}(\alpha)\Gamma_{q}(\gamma-\alpha)}\int_{0}^{1}t^{\alpha+n+k-1}\frac{(qt;q)_{\infty}}{(tq^{\gamma-\alpha};q)_{\infty}}d_{q}t,\label{2.57}
\end{split}
\end{eqnarray}
for $0<\Re(\alpha)<\Re(\gamma)$, $\gamma-\alpha\neq 0, -1, -2, -3, \ldots$, $Re(\alpha)>0$ and $n+k\geq 0$, we obtain an $q$-integral representation for $\mathbf{\Phi}_{1}$ (\ref{2.55}). Similarly, we obtain (\ref{2.56}).
\end{proof}
\begin{theorem}
For $r\in \mathbb{N}$, the differentiation formulas for $\mathbf{\Phi}_{1}$ $\mathbf{\Phi}_{2}$ and $\mathbf{\Phi}_{3}$ hold true
\begin{eqnarray}
\begin{split}
\mathbb{D}_{x,q}^{r}\bigg{[}x^{\beta+r-1}\mathbf{\Phi}_{1}(q^{\alpha},q^{\beta};q^{\gamma};q,x,y)\bigg{]}=\frac{(q^{\beta};q)_{r}}{(1-q)^{r}}x^{\beta-1}\mathbf{\Phi}_{1}(q^{\alpha},q^{\beta+r};q^{\gamma};q,x,y),\label{2.58}
\end{split}
\end{eqnarray}\vspace{-12pt}
\begin{eqnarray}
\begin{split}
\bigg{(}x^{2}\mathbb{D}_{x,q}\bigg{)}^{r}\bigg{[}x^{\alpha-r+1}\mathbf{\Phi}_{1}(q^{\alpha},q^{\beta};q^{\gamma};q,x,xy)\bigg{]}=\frac{(q^{\alpha};q)_{r}}{(1-q)^{r}}x^{\alpha+r}\mathbf{\Phi}_{1}(q^{\alpha+r},q^{\beta};q^{\gamma};q,x,xy),\\
\bigg{(}y^{2}\mathbb{D}_{y,q}\bigg{)}^{r}\bigg{[}y^{\alpha-r+1}\mathbf{\Phi}_{1}(q^{\alpha},q^{\beta};q^{\gamma};q,xy,y)\bigg{]}=\frac{(q^{\alpha};q)_{r}}{(1-q)^{r}}y^{\alpha+r}\mathbf{\Phi}_{1}(q^{\alpha+r},q^{\beta};q^{\gamma};q,xy,y),\label{2.59}
\end{split}
\end{eqnarray}\vspace{-12pt}
\begin{eqnarray}
\begin{split}
\mathbb{D}_{x,q}^{r}\bigg{[}x^{\alpha+r-1}\mathbf{\Phi}_{2}(q^{\alpha},q^{\beta};q^{\gamma};q,x,y)\bigg{]}=\frac{(q^{\alpha};q)_{r}}{(1-q)^{r}}x^{\alpha-1}\mathbf{\Phi}_{2}(q^{\alpha+r},q^{\beta};q^{\gamma};q,x,y),\\
\mathbb{D}_{y,q}^{r}\bigg{[}y^{\beta+r-1}\mathbf{\Phi}_{2}(q^{\alpha},q^{\beta};q^{\gamma};q,x,y)\bigg{]}=\frac{(q^{\beta};q)_{r}}{(1-q)^{r}}y^{\beta-1}\mathbf{\Phi}_{2}(q^{\alpha},q^{\beta+r};q^{\gamma};q,x,y)\label{2.60}
\end{split}
\end{eqnarray}
and
\begin{eqnarray}
\begin{split}
\mathbb{D}_{x,q}^{r}\bigg{[}x^{\alpha+r-1}\mathbf{\Phi}_{3}(q^{\alpha},q^{\beta};q^{\gamma};q,x,y)\bigg{]}=\frac{(q^{\alpha};q)_{r}}{(1-q)^{r}}x^{\alpha-1}\mathbf{\Phi}_{3}(q^{\alpha+r},q^{\beta};q^{\gamma};q,x,y).\label{2.61}
\end{split}
\end{eqnarray}
\end{theorem}
\begin{proof}
Using (see~\cite{gr1})
\begin{eqnarray*}
\begin{split}
&\mathbb{D}_{x,q}^{r}\bigg{[}x^{\beta+n+r-1}\bigg{]}=\frac{(q^{\beta+n};q)_{r}}{(1-q)^{r}}x^{\beta+n-1}
\end{split}
\end{eqnarray*}\vspace{-12pt}
and
\begin{eqnarray*}
\begin{split}
(q^{\beta};q)_{n}(q^{\beta+n};q)_{r}=(q^{\beta};q)_{n+r}=(q^{\beta};q)_{r}(q^{\beta+r};q)_{n},
\end{split}
\end{eqnarray*}
we get 
\begin{eqnarray*}
\begin{split}
&\mathbb{D}_{x,q}^{r}\bigg{[}x^{\beta+r-1}\mathbf{\Phi}_{1}(q^{\alpha},q^{\beta};q^{\gamma};q,x,y)\bigg{]}=\sum_{\ell,k=0}^{\infty}\frac{(q^{\alpha};q)_{n+k}(q^{\beta+n};q)_{r}(q^{\beta};q)_{n}}{(q^{\gamma};q)_{n+k}(q;q)_{n}(q;q)_{k})(1-q)^{r}}x^{\beta+n-1}y^{k}\\
&=\sum_{\ell,k=0}^{\infty}\frac{(q^{\alpha};q)_{n+k}(q^{\beta};q)_{n+r}}{(q^{\gamma};q)_{n+k}(q;q)_{n}(q;q)_{k})(1-q)^{r}}x^{\beta+n-1}y^{k}=\sum_{\ell,k=0}^{\infty}\frac{(q^{\alpha};q)_{n+k}(q^{\beta};q)_{r}(q^{\beta+r};q)_{n}}{(q^{\gamma};q)_{n+k}(q;q)_{n}(q;q)_{k})(1-q)^{r}}x^{\beta+n-1}y^{k}\\
&=\frac{(q^{\beta};q)_{r}}{(1-q)^{r}}\sum_{\ell,k=0}^{\infty}\frac{(q^{\alpha};q)_{n+k}(q^{\beta+r};q)_{n}}{(q^{\gamma};q)_{n+k}(q;q)_{k}(q;q)_{n})}x^{\beta+n-1}y^{k}.
\end{split}
\end{eqnarray*}
Proofs of Equations (\ref{2.59})--(\ref{2.61}) are along the same lines as of Equation (\ref{2.58}).
\end{proof}
\begin{remark}
The basic functions $\mathbf{\Phi}_{1}$, $\mathbf{\Phi}_{2}$ and $\mathbf{\Phi}_{3}$ are a $q$-analogue of the Humbert hypergeometric functions $\Phi_{1}$, $\Phi_{2}$ and $\Phi_{3}$ defined by (\cite{emot1}, page 225, Equations (20)--(22)):
\begin{eqnarray}
\begin{split}
\lim_{q\longrightarrow 1}\mathbf{\Phi}_{1}(q^{\alpha},q^{\beta};q^{\gamma};q,x,(1-q)y)=&\mathbf{\Phi}_{1}(\alpha,\beta;\gamma;x,y),\label{2.62}
\end{split}
\end{eqnarray}\vspace{-12pt}
\begin{eqnarray}
\begin{split}
&\lim_{q\longrightarrow 1}\mathbf{\Phi}_{2}(q^{\alpha},q^{\beta};q^{\gamma};q,(1-q)x,(1-q)y)&=\mathbf{\Phi}_{2}(\alpha,\beta;\gamma;x,y)\label{2.63}
\end{split}
\end{eqnarray}\vspace{-12pt}
and
\begin{eqnarray}
\begin{split}
&\lim_{q\longrightarrow 1}\mathbf{\Phi}_{3}(q^{\alpha};q^{\beta};q,(1-q)x,(1-q)^{2}y)&=\mathbf{\Phi}_{3}(\alpha;\beta;x,y).\label{2.64}
\end{split}
\end{eqnarray}
\end{remark}
\begin{proof}
Using the limit 
\begin{equation*}
\begin{split}
\lim_{q\longrightarrow 1}\frac{(q^{\alpha};q)_{n}}{(1-q)^{n}}&=(\alpha)_{n},
\end{split}
\end{equation*}
we obtain (\ref{2.62})--(\ref{2.64}).
\end{proof}
\begin{theorem}
The following relationship holds true between the basic functions $\mathbf{\Phi}_{1}$, $\mathbf{\Phi}_{2}$ and $\mathbf{\Phi}_{3}$:
\begin{eqnarray}
\begin{split}
&\lim_{\beta\longrightarrow \infty}\mathbf{\Phi}_{2}(q^{\alpha},q^{\beta};q^{\gamma};q,x,y)&=\mathbf{\Phi}_{3}(q^{\alpha};q^{\gamma};q,x,y),\label{2.65}
\end{split}
\end{eqnarray}\vspace{-12pt}
\begin{eqnarray}
\begin{split}
&\lim_{\beta\longrightarrow -\infty}\mathbf{\Phi}_{2}(q^{\alpha},q^{\beta};q^{\gamma};q,xq^{-\beta},y)&=\mathbf{\Phi}_{3}(q^{\alpha};q^{\gamma};q,-xq^{\frac{1}{2}(n-1)},y)\label{2.66}
\end{split}
\end{eqnarray}\vspace{-12pt}
and
\begin{eqnarray}
\begin{split}
\lim_{\alpha\longrightarrow \infty}\mathbf{\Phi}_{1}(q^{\alpha},q^{\beta};q^{\gamma};q,x,y)=&\mathbf{\Phi}_{3}(q^{\beta};q^{\gamma};x,y).\label{2.67}
\end{split}
\end{eqnarray}
\end{theorem}
\begin{proof}
Using the limit
\begin{equation*}
\begin{split}
\lim_{\alpha\longrightarrow \infty}(q^{\alpha};q)_{n}&=(0;q)_{n}=1,
\end{split}
\end{equation*}
we obtain (\ref{2.65})--(\ref{2.67}).
\end{proof}
\section{Concluding Remarks and Observations}
Due to the above investigation, our results are very remarkable and most general in their nature and are capable of giving many $q$-partial derivative formulas, $q$-differential formulas, $q$-recursion formulas, $q$-differential recursion relations, $q$-partial differential equations, summation formulas, transformation formulas, integral representations and interesting transformations results containing distinct kinds of $q$ or basic Humbert hypergeometric functions by some proper selections of symmetry parameters which are involved in the main results. Consequently, the results derived in this study possess some great applications in several different directions of mathematics, physics, statistics and engineering. Certain applications to other research topics, as well as further investigation of the properties and applications of these newly introduced extensions of basic Humbert hypergeometric functions, are left for future investigation by the authors and interested researchers.
\section*{Data Availability Statement} Data sharing is not applicable to this article as no data sets were generated or analyzed during the current study.
\subsection*{Conflicts of Interest} The author declare no conflicts of interest.
\section*{Author Contributions} The author has equal contributions. The author read and approved the final manuscript.
\section*{Funding} No funding was received for conducting this research.

\end{document}